\documentclass[12pt,oneside]{article}

\usepackage{amssymb,amsthm,amsmath}
\usepackage{pstricks,pst-node,pst-tree}
\usepackage{graphicx, pstcol, pst-plot}
\usepackage{subfigure}
\usepackage{color}
\definecolor{darkblue}{rgb}{0, 0, .4}
\definecolor{externblue}{rgb}{0, 0, .7}
\usepackage[bookmarks]{hyperref}
\hypersetup{
	colorlinks=true,
	linkcolor=darkblue,
	anchorcolor=darkblue,
	citecolor=darkblue,
	pagecolor=darkblue,
	urlcolor=externblue,
	pdftitle={The M\"obius function of the composition poset},
	pdfsubject={Combinatorics},
	pdfauthor={Bruce E. Sagan and Vincent Vatter},
	pdfkeywords={todo}
}
	
\newtheorem{theorem}{Theorem}[section]
\newtheorem{proposition}[theorem]{Proposition}
\newtheorem{lemma}[theorem]{Lemma}

\newtheorem{corollary}[theorem]{Corollary}

\newtheorem{conjecture}[theorem]{Conjecture}
\newtheorem{question}[theorem]{Question}

\newcommand{\ben}{\begin{enumerate}}
\newcommand{\een}{\end{enumerate}}
\newcommand{\ble}{\begin{lemma}}
\newcommand{\ele}{\end{lemma}}
\newcommand{\bth}{\begin{theorem}}
\renewcommand{\eth}{\end{theorem}}
\newcommand{\bpr}{\begin{proposition}}
\newcommand{\epr}{\end{proposition}}
\newcommand{\bco}{\begin{corollary}}
\newcommand{\eco}{\end{corollary}}
\newcommand{\bcon}{\begin{conj}}
\newcommand{\econ}{\end{conj}}
\newcommand{\bde}{\begin{defn}}
\newcommand{\ede}{\end{defn}}
\newcommand{\bex}{\begin{exa}}
\newcommand{\eex}{\end{exa}}
\newcommand{\barr}{\begin{array}}
\newcommand{\earr}{\end{array}}
\newcommand{\btab}{\begin{tabular}}
\newcommand{\etab}{\end{tabular}}
\newcommand{\beq}{\begin{equation}}
\newcommand{\eeq}{\end{equation}}
\newcommand{\bea}{\begin{eqnarray*}}
\newcommand{\eea}{\end{eqnarray*}}
\newcommand{\bce}{\begin{center}}
\newcommand{\ece}{\end{center}}
\newcommand{\bpi}{\begin{picture}}
\newcommand{\epi}{\end{picture}}
\newcommand{\bfi}{\begin{figure} \begin{center}}
\newcommand{\efi}{\end{center} \end{figure}}
\newcommand{\capt}{\caption}
\newcommand{\bsl}{\begin{slide}{}}
\newcommand{\esl}{\end{slide}}

\newcommand{\eqed}[1]{$\textcolor{white}{\qed}\hfill{#1}\hfill\qed$}

\newcommand{\ol}{\overline}

\newcommand{\hso}[1]{\hspace{-1pt}}
\newcommand{\vs}[1]{\vspace{#1}}

\newcommand{\emp}{\emptyset}

\newcommand{\sbs}{\subset}
\newcommand{\sbe}{\subseteq}

\newcommand{\zh}{\hat{0}}

\newcommand{\Ph}{\hat{P}}

\newcommand{\case}[4]
{\left\{\barr{ll}#1&\mbox{#2}\\#3&\mbox{#4}\earr\right.}

\def\<{\langle}
\def\>{\rangle}

\newcommand{\spn}[1]{\langle{#1}\rangle}

\newcommand{\ra}{\rightarrow}

\newcommand{\si}{\sigma}

\newcommand{\De}{\Delta}

\newcommand{\La}{\Lambda}

\newcommand{\bc}{{\bf c}}
\newcommand{\bd}{{\bf d}}

\newcommand{\bbP}{{\mathbb P}}

\newcommand{\bbZ}{{\mathbb Z}}

\newcommand{\cI}{{\mathcal I}}

\newcommand{\cJ}{{\mathcal J}}
\newcommand{\cK}{{\mathcal K}}

\newcommand{\bti}{\tilde{b}}

\newcommand{\mti}{\tilde{m}}

\newcommand{\chit}{\tilde{\chi}}

\newcommand{\sib}{\ol{\si}}
\newcommand{\taub}{\ol{\tau}}

\newcommand{\dil}{\displaystyle}

\newcounter{todocounter}
\newcommand{\todo}[1]{
	\addtocounter{todocounter}{1}
	\bigskip
	\noindent{\bf $\ll$ To-do \#\arabic{todocounter}:\rule{10pt}{0pt}#1
$\gg$}\bigskip
}

\newcommand{\vbar}{\overline{v}}
\newcommand{\vbarbar}{\overline{\overline{v}}}
\newcommand{\sbar}{\overline{s}}
\newcommand{\kbar}{\overline{k}}
\newcommand{\Supp}{\operatorname{Supp}}
\newcommand{\edge}[1]{\stackrel{#1}{\rule{20pt}{1pt}}}
\newcommand{\edg}{\edge{}}

\begin{document}

\title{The M\"obius function of the composition poset}

\author{Bruce E. Sagan\thanks{This work was partially done while the
author was on leave at DIMACS}\\[-5pt]
\small Department of Mathematics\\[-5pt]
\small Michigan State University\\[-5pt]
\small East Lansing, MI\\[-5pt]
\small \url{http://www.math.msu.edu/~sagan}\\[-5pt]
\small  
\href{mailto:sagan@math.msu.edu}{\texttt{sagan@math.msu.edu}}\\[6pt]
Vincent Vatter\thanks{Partially supported by an award from DIMACS and  
an NSF VIGRE grant to the Rutgers University Department of  
Mathematics.}\\[-5pt]
\small Department of Mathematics\\[-5pt]
\small Rutgers University\\[-5pt]
\small Piscataway, NJ\\[-5pt]
\small \url{http://math.rutgers.edu/~vatter/}\\[-5pt]
\small  
\href{mailto: 
vatter@math.rutgers.edu}{\texttt{vatter@math.rutgers.edu}}}

\date{\today \\[6pt]
	\begin{flushleft}
	\small Key Words: composition, discrete Morse function, M\"obius  
function,
permutation pattern, subword order\\[6pt]
	\small AMS classifications: 06A07, 05E25, 68R15
	\end{flushleft}
           }

\maketitle

\begin{abstract}
We determine the M\"obius function of the poset of compositions of an
integer.  In fact we give two proofs of this formula, one using an
involution and one involving discrete Morse theory.  The composition
poset turns out to be intimately connected with subword order, whose
M\"obius function was determined by Bj\"orner.  We show that
using a generalization of subword order, we can obtain both Bj\"orner's
results and our own as special cases.
\end{abstract}

If $A$ is any set then the corresponding {\it Kleene closure\/} or
{\it  free monoid}, $A^*$ is the set of words with letters from $A$,  
i.e.,
$$
A^*=\{w=w(1)w(2)\ldots w(n)\ |\ \mbox{$n\ge0$ and $w(i)\in A$ for all
$i$}\}.
$$
We denote the length (number of elements) of $w$ by $|w|$.

Letting $\bbP$ denote the positive integers, we see that $\bbP^*$ is
the set of integer compositions (ordered partitions).  We can
turn $\bbP^*$ into a partially ordered set by letting $u\le w$ if
there is a subword $w(i_1) w(i_2)\ldots w(i_l)$ of $w$ having length
$l=|u|$ such that $u(j)\leq w(i_j)$ for $1\leq j\leq l$.  For example,
$433 \leq 16243$ because of the subword $643$.
Bergeron, Bousquet-M\'elou and Dulucq~\cite{bbd:spc} were the first to
study
$\bbP^*$, enumerating saturated chains that begin at its minimal
element.  Snellman has also studied saturated chains in this poset as
well as two other partial orders on
$\bbP^*$~\cite{sne:scc,sne:spa}.  One of the main
results of this paper is a formula for the M\"obius function of
$\bbP^*$.

This order on $\bbP^*$ is closely related to subword order.  If
$A$ is any set then the {\it subword order\/} on $A^*$ is defined by
letting $u\le w$ if $w$ contains a subsequence
$w(i_1), w(i_2),\ldots,  w(i_l)$
such that $u(j)= w(i_j)$ for $1\leq j\leq l=|u|$.  By way of
illustration,
if $A=\{a,b\}$ then $abba\le aabbbaba$ because
$w(1)w(3)w(5)w(8)=abba$.  Note that we use the
notation $A^*$ when referring to subword order as opposed to the
partial order
on $\bbP^*$, even though we use $\le$ for both.  We will always
give enough context to make it clear which poset we are dealing with.
Bj\"orner~\cite{bjo:mfs} was the first to completely determine the
M\"obius function of subword order, although special cases had been
obtained previously by Farmer~\cite{far:chp} and
Viennot~\cite{vie:mcs}. In fact, Bj\"orner gave two proofs of his
formula, one using an involution~\cite{bjo:mff} and one using
shellability~\cite{bjo:mfs}.  He also gave a demonstration with
Reutenauer~\cite{br:rmf} via  generating functions on monoids.
Another proof  was given by Warnke~\cite{war:mfs} using induction
while Wang and Ma~\cite{wm:gcp} used Cohen-Macaulayness to investigate  
$\mu$.
We derive the M\"obius function formula for $\bbP^*$ using first
combinatorial and then topological techniques.

The rest of this paper is structured as follows.  In the next section
we review Bj\"orner's result for subword order as well as the
related definitions which will be useful for $\bbP^*$.
It contains a statement of our formula for the M\"obius function of
$\bbP^*$ in Theorem~\ref{mobius-main}.  Section~\ref{psi}
is devoted to giving a proof of this theorem using a
sign-reversing involution.

Although intervals in $A^*$ are shellable,
those in $\bbP^*$ need not even be connected as evident from the
example in Figure~\ref{322}, so we need a more powerful tool to study
the topology of the composition poset. For this we turn to discrete
Morse
theory, which was developed by Forman~\cite{for:dmt,for:mtc} and
can be used to compute the homology of any CW-complex. A method
for applying this theory to the order complex of a poset was given by
Babson and Hersh~\cite{bh:dmf} and further studied by Hersh
herself~\cite{her:odf}.  Since this is a relatively
recent addition to the combinatorial toolbox, we provide an
exposition of the basic ideas of the theory in Section~\ref{idm}.  The
subsequent section gives a Morse theoretic proof of
Theorem~\ref{mobius-main}.

The similarity between the formulas for the M\"obius functions of
$A^*$ and $\bbP^*$ leads one to ask if there is a common
generalization.  In fact, if $P$ is any poset then there is a partial
order on $P^*$ which we call generalized subword order.  It has been
used in the context of well-quasi-ordering;
see Kruskal's article~\cite{kru:twq} for a survey of the early  
literature.
When $P$ is a chain or an antichain in the present context,
one recovers our results or Bj\"orner's, respectively.
This construction is studied
in Section~\ref{gso}.  Finally, we end with a section of comments and
open problems.

\section{Subword and composition order}
\label{sco}

We now review Bj\"orner's formula for the M\"obius function of
subword order, reformulating it slightly so as to emphasize the
connection to the composition order which is our main objective.
We assume the reader is familiar with M\"obius functions, but all
the necessary definitions and theorems we use here can be found
in Stanley's text~\cite[\S 3.6--3.7]{sta:ec1}.

We first need to restate the definition of the partial order in
$A^*$ in a way that, although slightly more complicated, has a
direct connection with the M\"obius function.
Suppose we have a distinguished symbol, $0$, and suppose that
$0\not\in A$.  Then
a word $\eta=\eta(1)\eta(2)\ldots\eta(n)\in (A\cup 0)^*$ has {\it
support\/}
$$
\Supp \eta =\{i\ |\ \eta(i)\neq 0\}.
$$
An {\it expansion\/} of $u\in (A\cup 0)^*$ is a word $\eta_u\in(A\cup  
0)^*$ such that the restrictions of $u$ and $\eta_u$ to their supports
are equal.  For example, if $u=abba$
then one expansion of $u$ is $\eta_u=a0b0b00a$.  An {\it
embedding\/} of $u$ into $w$ is an expansion $\eta_{uw}$ of $u$ which
has length $|w|$ and satisfies
$$
\mbox{$\eta_{uw}(i)=w(i)$ for all $i\in\Supp\eta_{uw}$.}
$$
Note that $u\le w$ in $A^*$ if and only if there is an embedding of
$u$ into $w$.  The example expansion $\eta_{uw}$ given above is exactly
the
embedding which corresponds to the subword of
$w=aabbbaba$ given at the beginning of the third
paragraph of this paper.  If $w$ is clear from context we simply
write $\eta_u$ for $\eta_{uw}$.

The M\"obius function of $A^*$ counts certain types of
embeddings.  If $a\in A$ then a {\it run\/} of $a$'s in $w$ is a
maximal interval of indices $[r,t]$ such that
$$
w(r)=w(r+1)=\ldots=w(t)=a.
$$
Continuing our example, the runs in $w=aabbbaba$
are $[1,2]$, $[3,5]$, $[6,6]$, $[7,7]$, and $[8,8]$.  Call an
embedding $\eta_u$ into $w$ {\it normal\/} if, for every $a\in A$ and
every run $[r,t]$ of $a$'s, we have
$$
(r,t] \sbe \Supp \eta_u
$$
where $(r,t]$ denotes the half-open interval.  In our running
example, this means that the second $a$ and the fourth and fifth $b$'s
must be in any normal embedding.  (If $r=t$ then $(r,t]=\emp$, so
there is no restriction on runs of one element.)  Thus in this case
there are precisely two normal embeddings of $u$ into $w$, namely
$$
\mbox{$\eta_u = 0a0bba00$ and $0a0bb00a$.}
$$
Define $\binom{w}{u}_n$ to be the number of normal embeddings of $u$  
into $w$.

We now have everything in place to state Bj\"orner's result.
\bth[Bj\"orner~\cite{bjo:mfs}]
\label{bjo}
If $u,w\in A^*$ then\\
\eqed{\mu(u,w)=(-1)^{|w|-|u|}\dil{w\choose u}_n.}
\eth

Finishing our example, we see that
$$
\mu(abba,aabbbaba)=(-1)^{8-4} \cdot 2 = 2.
$$

We now turn to $\bbP^*$.  The definitions of support and expansion are
exactly as before, but the notion of embedding must be updated
to reflect the different partial order.
To this end we define an {\it embedding\/} of $u$ into
$w$ as an expansion $\eta_u$ of $u$ having length $|w|$ such that
$$
\mbox{$\eta_u(i)\le w(i)$ for $1\le i\le |w|$.}
$$
Again, $u\le w$ in $\bbP^*$ is equivalent to the existence of
an embedding
of $u$ into $w$.  An interval in $\bbP^*$ is displayed in
Figure~\ref{322}.

If $u\le w$ and $\eta_w$ is an expansion of $w$ then there is a unique
{\it last\/} or {\it rightmost\/} embedding $\rho_{uw}$ of $u$ into
$\eta_w$ which has the property that for any other embedding
$\eta_{uw}$ of $u$ into $\eta_w$ one has $\Supp(\eta_{uw})\le
\Supp(\rho_{uw})$.  (If $S=\{i_1<\cdots<i_m\}$ and
$S'=\{i_1'<\cdots<i_m'\}$ then we write $S\le S'$ to mean that
$i_j\le i_j'$ for $1\le j \le m$.)  Note that $\rho_{uw}$ depends on
$\eta_w$,
not just on $w$, but the expansion of $w$ used will always be clear
from context.

Like subword order, we define normal embeddings for $\bbP^*$ in
terms of runs (defined in the same way in this context).  We
say that an embedding $\eta_u$ into $w$ is {\it normal\/} if the
following conditions hold.
\begin{enumerate}
\item\label{normal1} For $1\le i\le |w|$ we have $\eta_u(i)=w(i)$,
$w(i)-1$, or $0$.\vs{3pt}
\item\label{normal2} For all $k\ge 1$ and every run $[r,t]$ of $k$'s
in $w$, we have\vs{3pt}
\begin{enumerate}
\item $(r,t]\sbe \Supp \eta_u$ if  $k=1$,\vs{3pt}
\item $r\in\Supp \eta_u$ if $k\ge 2$.
\end{enumerate}
\end{enumerate}
Comparing this with Bj\"orner's definition, we see that in $\bbP^*$ a
normal
embedding can have three possible values at each
position instead of two.  Also, the run condition for ones is the same
as in $A^*$, while that condition for integers greater than one is
complementary.  As an example, if $w=2211133$ and
$u=21113$, then there are two normal embeddings, namely
$\eta_u= 2101130$ and $2011130$.  Note that
$2001113$  and $0211130$ are not normal since
they violate conditions~\eqref{normal1} and~\eqref{normal2},  
respectively.

The sign of a normal embedding depends on the embedding itself and
not just the length of the compositions.  Given a normal embedding
$\eta_u$ into $w$ we define its {\it defect} to be
$$
d(\eta_u)=\#\{i\ |\ \eta_u(i)=w(i)-1\}.
$$

We can now state our main theorem about $\bbP^*$.
\begin{theorem}
\label{mobius-main}
If $u,w\in\bbP^*$ then
$$
\mu(u,w)=\sum_{\eta_u} (-1)^{d(\eta_u)}
$$
where the sum is over all normal embeddings $\eta_u$ into $w$.
\end{theorem}

In the example of the previous paragraph, this gives
$$
\mu(21113,2211133) = (-1)^2+(-1)^0=2.
$$
Although this example does not show it, it is possible to have
cancellation among the terms in the sum for $\mu$.

\section{Proof by sign-reversing involution}
\label{psi}

We now prove Theorem~\ref{mobius-main} using a sign-reversing
involution.  The proof is similar in nature to Bj\"orner's proof
in~\cite{bjo:mff}, but is significantly more complicated.

\medskip

\noindent {\it Proof (of Theorem~\ref{mobius-main}).}
If $u=w$ then there is exactly one normal embedding $\eta_{ww}$ and it
has defect $0$.  This gives $(-1)^0=1=\mu(w,w)$, as desired.

Now assume that $u<w=w(1)\cdots w(n)$.  Since the M\"obius recurrence
uniquely defines $\mu$, it suffices to show that
$$
\sum_{v\in[u,w]}\sum_{\eta_{vw}} (-1)^{\eta_{vw}}=0,
$$
where the inner sum is over all normal embeddings of $v$ into $w$.
We prove this by constructing a sign-reversing involution on the set of
normal embeddings $\eta_{vw}$ for $v\in[u,w]$.

Let $\rho_{uv}$ denote the rightmost embedding of $u$ into $\eta_{vw}$,
so for all $i\in[1,n]$,
\begin{equation}\label{inv-rightmost}
\rho_{uv}(i)\le \eta_{vw}(i)\le w(i).
\end{equation}
Also let $s$ denote the left-most index where $\rho_{uv}$ and $w$
differ, i.e., $s=\min\{i : \rho_{uv}(i)\neq w(i)\}$.  Since $u<w$, $s$
must exist, and by the definition of $s$ we have
\begin{equation}\label{def-s-consequence}
\rho_{uv}(i)=\eta_{vw}(i)=w(i)\mbox{ for all $i\in[1,s)$.}
\end{equation}
Set $k=w(s)$, so $\eta_{vw}(s)$ is $k$, $k-1$, or $0$ by normality and
$\rho_{uv}(s)\le k-1$.  Finally, let $[r,t]$ denote the indices of the
run of $k$'s in $w$ that contains the index $s$.

Our involution maps $\eta_{vw}$ to the embedding $\eta_{\vbar
    w}=\eta_{\vbar w}(1)\cdots\eta_{\vbar w}(n)$ where $\eta_{\vbar
    w}(i)=\eta_{vw}(i)$ for all $i\neq s$ and $\eta_{\vbar w}(s)$ is
determined by the following rules.  If $k=1$ then
\begin{equation}\label{inv1}
\eta_{\vbar w}(s)
=
\left\{
\begin{array}{ll}
0&\mbox{if $\eta_{vw}(s)=1$,}\\
1&\mbox{if $\eta_{vw}(s)=0$.}
\end{array}
\right.
\end{equation}
If $k\ge 2$, $s>r$, and $\rho_{uv}(s)=0$ then
\begin{equation}\label{inv2}
\eta_{\vbar w}(s)
=
\left\{
\begin{array}{ll}
0&\mbox{if $\eta_{vw}(s)=k-1$,}\\
k-1&\mbox{if $\eta_{vw}(s)=0$.}\\
\end{array}
\right.
\end{equation}
Finally, if $k\ge 2$ and either $s=r$ or $\rho_{uv}(s)\neq 0$ then
\begin{equation}\label{inv3}
\eta_{\vbar w}(s)
=
\left\{
\begin{array}{ll}
k-1&\mbox{if $\eta_{vw}(s)=k$,}\\
k&\mbox{if $\eta_{vw}(s)=k-1$.}
\end{array}
\right.
\end{equation}

It is not obvious that this map is
defined for all normal embeddings $\eta_{vw}$ when $k\ge 2$.  For
example, if $s>r$ and $\rho_{uv}(s)=0$, then we should apply
\eqref{inv2}, but it is a priori possible that $\eta_{vw}(s)=k$, in
which case \eqref{inv2} is not defined.  However, since $s>r$,
$\rho_{uv}(s-1)=w(s-1)=k$ by \eqref{def-s-consequence}, and this
contradicts our choice of $\rho_{uv}$ as the rightmost embedding of
$u$ into $\eta_{vw}$.  Similar issues arise when $s=r$ or
$\rho_{uv}(s)\neq 0$: if $s=r$ then $\eta_{vw}(s)\neq 0$ by normality,
and if $\rho_{uv}(s)\neq 0$ then $\eta_{vw}(s)\neq 0$ by
\eqref{inv-rightmost}.

Having established that this map is indeed defined on all normal
embeddings, we have several properties to prove.  First, it is evident
from \eqref{inv1}, \eqref{inv2}, and \eqref{inv3} that the number of
elements equal to $k-1$ changes by exactly one in passing from
$\eta_{vw}$ to $\eta_{\vbar w}$, and thus $(-1)^{\eta_{\vbar
      w}}=-(-1)^{\eta_{vw}}$.  We must now prove that $\eta_{\vbar w}$
is a normal embedding of $\vbar$ into $w$ for some
$\vbar\in[u,w]$ and that this map is an involution.

We begin by showing that $\eta_{\vbar w}$ is an embedding of some
$\vbar\in[u,w]$ into $w$.  It follows from the fact that $\eta_{vw}$ is
an embedding into $w$ and the definition of our map that $\eta_{\vbar
    w}$ is an embedding of some word
$\vbar$ into $w$, so we only need to show that $\vbar\ge u$.  We
prove this by showing that
\begin{equation}\label{vgeu}
\eta_{\vbar w}(i)\ge\rho_{uv}(i)
\end{equation}
for all $i\in[1,n]$.  This is clear for all $i\neq s$ because for these
indices $\eta_{\vbar w}(i)=\eta_{vw}(i)\ge\rho_{uv}(i)$.  Furthermore,
$\rho_{uv}(s)\le k-1$ by the definition of $s$, so the only case in
which \eqref{vgeu} is not immediate is when $k\ge 2$ and $\eta_{\vbar
    w}(s)=0$.  However, this can only occur from using \eqref{inv2},  
which
requires that $\rho_{uv}(s)=0$, completing the demonstration that
$\vbar\in[u,w]$.

We  now aim to show that $\eta_{\vbar w}$ is a normal embedding.  If
$\eta_{\vbar w}(s)>0$ then $\Supp(\eta_{\vbar w})\supseteq
\Supp(\eta_{vw})$, so the normality of $\eta_{\vbar w}$ follows from
the normality of $\eta_{vw}$ and the fact that $\eta_{\vbar w}(s)$ is
either $k-1$ or $k$.  If $\eta_{\vbar w}(s)=0$ then there are two
cases depending on whether \eqref{inv1} or \eqref{inv2} was applied.
Suppose first that \eqref{inv1} was applied, so $k=1$.  Comparing
\eqref{inv1} with the definition of normality, we see that it suffices
to show $s=r$.  Suppose to the contrary that $s\in(r,t]$.  Since $s$
is the left-most position at which $\rho_{uv}$ and $w$ differ, we have
$\rho_{uv}(s)=0$ and $\rho_{uv}(s-1)=1$.  However, this contradicts
our choice of $\rho_{uv}$ as the rightmost embedding of $u$ in
$\eta_{vw}$.  Now suppose that \eqref{inv2} was applied, so $s>r$ and
$\rho_{uv}(s)=0$.  Then \eqref{def-s-consequence} implies that
$\eta_{\vbar w}(r)=\eta_{vw}(r)=k$, and normality is preserved.

It only remains to show that this map is an involution.  Consider
applying the map to $\eta_{\vbar w}$.  In this process we define
$\rho_{u\vbar}$ to be the rightmost embedding of $u$ into $\eta_{\vbar
    w}$, $\sbar=\min\{i : \rho_{u\vbar}(i)\neq w(i)\}$, and
$\kbar=w(\sbar)$.  We then follow the rules \eqref{inv1},
\eqref{inv2}, and \eqref{inv3} to construct an embedding
$\eta_{\vbarbar w}$, which we would like to show is equal to
$\eta_{vw}$.

First we claim that $\rho_{u\vbar}=\rho_{uv}$.  Suppose to the
contrary that $\rho_{u\vbar}\neq\rho_{uv}$.  By \eqref{vgeu},
$\eta_{\vbar w}(i)\ge\rho_{uv}(i)$ for all $i$, so $\rho_{uv}$ also
gives an embedding of $u$ into $\vbar$, and thus the only way we can
have $\rho_{u\vbar}\neq\rho_{uv}$ is if $\rho_{u\vbar}$ is further to
the right than $\rho_{uv}$.  This requires that
\beq
\label{r}
0=\rho_{uv}(s)<\rho_{u\vbar}(s)\le k
\eeq
and that
\beq
\label{e}
\eta_{vw}(s)<\eta_{\vbar w}(s).
\eeq

Because we are assuming that $\rho_{\vbar w}$ is further to the right
than $\rho_{vw}$, there is some position to the left of $s$ at which
$\rho_{uv}$ is nonzero.  In fact, \eqref{def-s-consequence} shows that
we must have $\rho_{uv}(s-1)\neq 0$ and also implies that
\beq
\label{re}
\eta_{vw}(s)<\rho_{u\vbar}(s)=\rho_{uv}(s-1)=w(s-1).
\eeq
We now consider the three cases
arising from each of the rules \eqref{inv1}, \eqref{inv2}, and
\eqref{inv3} in turn.

Suppose \eqref{inv1} was applied so $k=1$.  It follows from~\eqref{e}
and the definition of our map that $\eta_{vw}(s)=0$.
Also~\eqref{r} and~\eqref{re} give
$w(s-1)=\rho_{u\vbar}(s)=1$.  But then $\eta_{vw}$ zeroed out a $1$
which was not the first in its run, contradicting normality.

Now suppose \eqref{inv2} was applied.  Then by~\eqref{e} we
have $\eta_{\vbar w}(s)=k-1$.
We also have $s>r$ which in conjunction with~\eqref{re}
gives $\rho_{u\vbar}(s)=w(s-1)=k$.  This implies that
$\eta_{\vbar w}(s)<\rho_{u\vbar}(s)$, but that contradicts
the $\vbar$ version of~\eqref{inv-rightmost}.

Finally suppose that \eqref{inv3} was used.  By~\eqref{r}
it must be the case that $s=r$.  Also,  equation~\eqref{e} gives
$\eta_{vw}(s)=k-1$ and $\eta_{\vbar w}(s)=k$.  Now
applying~\eqref{r} and~\eqref{re} we have
$k-1<w(s-1)=\rho_{u\vbar}(s)\le k$, so $w(s-1)=k$ which
contradicts that fact that $s=r$.

Now that we have established the equality of $\rho_{uv}$ and
$\rho_{u\vbar}$, the fact that this map is an involution can be
readily observed.  We must have $\sbar=s$, so $\kbar=k$, and thus we
apply the same rule to go from $\eta_{\vbar w}$ to
$\eta_{\vbarbar w}$ as we applied to get $\eta_{\vbar w}$ from
$\eta_{\vbar w}$, and each of these rules is clearly an involution. \qed

\section{Introduction to discrete Morse theory}
\label{idm}
In this section we review the basic ideas behind Forman's discrete
Morse theory~\cite{for:dmt,for:mtc} as well as Babson and Hersh's
method for applying the theory to the order complex of a
poset~\cite{bh:dmf}.

Let $X$ be a CW-complex. Since we will be working in reduced
homology, we assume that  $X$  has an empty cell $\emp$ of
dimension $-1$ which is contained in every cell of $X$.  If $\si$ is a
$d$-cell (cell of dimension $d$) in $X$ then let $\si^\partial$ be the
set of $(d-1)$-cells $\tau$ which are contained in the closure $\sib$.
Dually, let $\si^\delta$ denote the set of $(d+1)$-cells $\tau$ such  
that
$\si\sbs\taub$.

A real-valued function $f$ on the cells of $X$ is a {\it Morse  
function\/}
if it satisfies the following two conditions.
\ben
\item  For every cell $\si$ of $X$ we have\vs{3pt}
   \ben
   \item $\#\{\tau\in\si^\partial\ |\ f(\tau)\ge f(\si)\}\le 1$,  
and\vs{3pt}
   \item $\#\{\tau\in\si^\delta\ |\ f(\tau)\le f(\si)\}\le 1$.\vs{3pt}
   \een
\item If $\tau\in\si^\delta$ and $f(\tau)\le f(\si)$ then $\si$ is a
   regular face of $\tau$.
\een
Intuitively the first condition says that, with only certain exceptions,
$f$ increases with dimension.  In fact, $f(\si)=\dim\si$ is a
perfectly good Morse function on $X$, although we will see shortly
that it is not very interesting.  A simple example of a Morse function
on a CW-complex is given in Figure~\ref{mf} where the value of $f$ is
given next to each cell $\si$ and we also set $f(\emp)=-1$.

\thicklines
\setlength{\unitlength}{2pt}
\bfi
\bpi(60,60)(0,0)
\pscircle(2.1,2.1){1.4}
\put(30,10){\circle*{3}}
\put(30,50){\circle*{3}}
\put(30,5){\makebox(0,0){$0$}}
\put(5,30){\makebox(0,0){$1$}}
\put(30,55){\makebox(0,0){$2$}}
\put(55,30){\makebox(0,0){$3$}}
\epi
\capt{A Morse function on a CW-complex}
\label{mf}
\efi

The fact that condition (1) holds for every cell implies that, in fact,  
at most one
of the two sets under consideration has cardinality equal to 1.  Thus
the function $f$ induces a {\it Morse matching\/} between pairs of
cells $\si,\tau$ with $\sigma\in\tau^\partial$ and $f(\sigma)\ge  
f(\tau)$.  The
regularity condition ensures that for each such pair there is an
elementary collapse of $\taub$ onto $\taub-(\tau\cup\si)$.
The cells which are not matched by $f$ are called {\it critical}.
Since each collapse is a homotopy equivalence,  $X$ can be collapsed  
onto
a homotopic complex $X^f$ built from the critical cells.  In our
example, the cells labeled $1$ and $2$ are matched and after
collapsing we clearly have a complex which is still homotopically a
circle.  Note that if we take $f$ to be the dimension function then
every cell is critical and $X^f=X$, so the cell complex does not
simplify in this case which does not help in understanding its  
structure.

Let $\mti_d$ be the number of critical $d$-cells of $X$ and let $\bti_d$
be the $d$-th reduced Betti number over the integers.  We also use
$\chit(X)$ for the reduced Euler characteristic.  From the
considerations in the previous paragraph, we have the following
{\it Morse inequalities\/} which are analogous to those in traditional
Morse theory.
\bth[Forman~\cite{for:mtc}]
\label{me}
For any Morse function on a cell complex $X$ we have
\ben
\item $\bti_d\le \mti_d$ for $d\ge-1$, and\vs{3pt}
\item $\chit(X)=\dil\sum_{d\ge -1}(-1)^d \mti_d$.\qed
\een
\eth
(One can get further inequalities relating various partial alternating
sums of the $\bti_d$ and $\mti_d$.)
Continuing our example, we see that $\mti_{-1}=\mti_0=\mti_1=1$ which  
bound
$\bti_{-1}=\bti_0=0$ and $\bti_1=1$, as well as $\chit(X)=-1+1-1=-1$, as
expected.

We now turn to the special case of order complexes.  Let $P$ be a
poset and consider an open interval $(u,w)$ in $P$.  The corresponding  
{\it
order complex\/} $\De(u,w)$ is the abstract simplicial complex
whose simplices (faces) are the chains in
$(u,w)$.  We are interested in the order complex because of the
fundamental fact~\cite{rot:tmf} that
\beq
\label{chit}
\mu(u,w)=\chit(\De(u,w)).
\eeq

Therefore finding a Morse function for $\De(u,w)$ could permit us to  
derive
the corresponding M\"obius value as well as give extra information
about its Betti numbers.  Suppose we have an ordering of the maximal
chains of $(u,w)$ (facets of $\De(u,w)$), say $C_1,C_2,\ldots,C_l$.
Call a face (subchain) $\si$ of $C_k$ {\it new\/} if it is not
contained in any $C_j$ for $j<k$.  We would like to construct a Morse
matching inductively,
where at the $k$th stage we extend the matching on the faces in $C_j$
for $j<k$ by matching up as many of the new faces $\si$ in $C_k$ as  
possible.
It turns out that under fairly mild conditions on the
facet ordering, one can construct such a matching so that all the new
faces in $C_k$ are matched if there are an even number of them, and  
only one is
left unmatched if the number is odd.  Thus adding each facet
contributes at most one critical cell.  A maximal chain contributing a
critical cell is called a {\it critical chain}.  In reading the details  
of this
construction, the reader may find it useful to refer to the example
of the interval $[322,3322]\sbs\bbP^*$ given in Figure~\ref{322}.
Note that by abuse of notation we include $u$ and $w$ when writing out
a maximal chain $C$, even though $C$ is really a subset of the open
interval $(u,w)$.
Also, because of the way our chain order is constructed, we
start with the top element $w$ and work down to $u$ which is dual to
what is done normally.
Thus in a chain $C$, terms like ``first'' and ``last'' refer to this
ordering of $C$'s  elements.
Finally, we list the elements of a chain as
embeddings into $w$ for reasons which will become apparent when we
also describe the labels given to the edges (covers) of a chain.

\begin{figure}
\begin{center}
\begin{footnotesize}
$$
\psset{nodesep=2pt,colsep=18pt,rowsep=25pt}
\begin{psmatrix}
\psspan{2}&[name=3322] 3322&\psspan{2}\\[0pt]
[name=2322] 2322&[name=3222] 3222&&
[name=3312] 3312&[name=3321] 3321\\[0pt]
[name=1322] 1322&[name=3122] 3122&
[name=3212] 3212&[name=3221] 3221&[name=332] 332\\[0pt]
\psspan{2}&[name=322] 322&\psspan{2}
\ncline{3322}{2322}
\ncline{3322}{3222}
\ncline{3322}{3312}
\ncline{3322}{3321}
\ncline{2322}{1322}
\ncline{3222}{3122}
\ncline{3222}{3212}
\ncline{3222}{3221}
\ncline{3321}{3221}
\ncline{3321}{332}
\ncline{3312}{3212}
\ncline{3312}{332}
\ncline{1322}{322}
\ncline{3122}{322}
\ncline{3212}{322}
\ncline{3221}{322}
\ncline{3212}{322}
\ncline{332}{322}
\end{psmatrix}
$$
\end{footnotesize}
\begin{footnotesize}
$$
\begin{psmatrix}[rowsep=3pt,colsep=5pt]
C_{1}:&\rule{10pt}{0pt}&3322&\edge{1}&2322&\edge{1}&1322&\edge{1}&0322\\
C_{2}:&\rule{10pt}{0pt}&3322&\edge{2}&3222&\edge{2}&3122&\edge{2}&3022\\
C_{3}:&\rule{10pt}{0pt}&3322&\edge{2}&3222&\edge{3}&3212&\edge{3}&3202\\
C_{4}:&\rule{10pt}{0pt}&3322&\edge{2}&3222&\edge{4}&3221&\edge{4}&3220\\
C_{5}:&\rule{10pt}{0pt}&3322&\edge{3}&3312&\edge{2}&3212&\edge{3}&3202\\
C_{6}:&\rule{10pt}{0pt}&3322&\edge{3}&3312&\edge{3}&3302&\edge{2}&3202\\
C_{7}:&\rule{10pt}{0pt}&3322&\edge{4}&3321&\edge{2}&3221&\edge{4}&3220\\
C_{8}:&\rule{10pt}{0pt}&3322&\edge{4}&3321&\edge{4}&3320&\edge{2}&3220\\
\ncbox[nodesep=0.07,boxsize=0.18,linearc=0.15,linestyle=solid,
linecolor=gray]{2,5}{2,7}
\ncbox[nodesep=0.07,boxsize=0.18,linearc=0.15,linestyle=solid,
linecolor=gray]{3,7}{3,7}
\ncbox[nodesep=0.07,boxsize=0.18,linearc=0.15,linestyle=solid,
linecolor=gray]{4,7}{4,7}
\ncbox[nodesep=0.07,boxsize=0.18,linearc=0.15,linestyle=solid,
linecolor=gray]{5,5}{5,5}
\ncbox[nodesep=0.07,boxsize=0.18,linearc=0.15,linestyle=solid,
linecolor=gray]{6,7}{6,7}
\ncbox[nodesep=0.07,boxsize=0.18,linearc=0.15,linestyle=solid,
linecolor=gray]{7,7}{7,7}
\ncbox[nodesep=0.07,boxsize=0.18,linearc=0.15,linestyle=solid,
linecolor=gray]{8,5}{8,5}
\ncbox[nodesep=0.07,boxsize=0.18,linearc=0.15,linestyle=solid,
linecolor=gray]{8,7}{8,7}
\end{psmatrix}
$$
\end{footnotesize}
\end{center}
\caption{The interval $[322,3322]$ and its maximal chains}
\label{322}
\end{figure}

To define the types of chain orderings we consider, suppose we
have two chains $C:\ w=v_0 \edg v_1 \edg \ldots \edg u$
and $C':\ w=v_0' \edg v_1' \edg \ldots \edg u$ where $x\edg y$
means that $x$ covers $y$.  Then we say that $C$ and $C'$ {\it agree
to index $j$\/} if $v_i=v_i'$ for $i\le j$.  In addition, $C$ and $C'$
{\it diverge from index $j$\/} if they agree to index $j$ and
$v_{j+1}\neq v_{j+1}'$.  In addition, we use the notation
$C'<C$ to mean that $C'$ comes before $C$ in the order under  
consideration.
An ordering of the maximal chains of $[u,w]$
is a {\it poset lexicographic order}, or {\it PL-order\/}  for short,
if it satisfies the following condition.
Suppose $C'$ and $C$ diverge from index $j$ with $C'<C$.
Then for any maximal chains $D'$ and $D$ which agree to
index $j+1$ with $C'$ and $C$, respectively, we must have $D'<D$.
Note that orderings coming from
the EL-labelings introduced by Bj\"orner~\cite{bjo:scp} or from the
more general CL-labelings of Bj\"orner and Wachs~\cite{bw:boc} are
PL-orders as long as one breaks ties among labels consistently.

The PL-order we use in $\bbP^*$ is as follows.  If $x\edg y$ is a
cover then $y$ is obtained from $x$ by reducing a single part of $x$
by $1$.  Thus there is a unique normal embedding of $y$ into $x$, since  
if
a $1$ is reduced to $0$ then it must be the first element in the run
of ones to which it belongs.  Similarly, for any expansion $\eta_x$
there is a unique normal embedding of $y$ into $\eta_x$.  Now given
any chain $C:\ w=v_0 \edg v_1 \edg \ldots \edg u$ we inductively
associate with each $v_j$ an embedding $\eta_{v_j}$ into $w$ where
$\eta_{v_0}=w$ and, for $j\ge0$, $\eta_{v_{j+1}}$ is the unique normal
embedding of $v_{j+1}$ into $\eta_{v_j}$.  We label the edge
$v_j\edg v_{j+1}$ of $C$ with the index $i$ of the position which was
decreased in passing from $\eta_{v_j}$ to $\eta_{v_{j+1}}$.
Furthermore, we often write $\eta_{v_j}$ in place of $v_j$ when
listing the elements of $C$.  Figure~\ref{322} illustrates this
labeling.
It is important to note that although $\eta_{v_{j+1}}$ is
normal in $\eta_{v_j}$, it need not be normal in $w$.  We should also
remark that this labeling is similar to the one used by
Bj\"orner~\cite{bjo:mfs} in his CL-shelling of the intervals in subword
order.  Finally, if one orders the chains of $[u,w]$
using ordinary lexicographic order on their label sequences, then the
result is a PL-order.  This is due to the fact that if two chains
agree to index $j$ then their first $j$ labels are the same.  The
chains in Figure~\ref{322} are listed in PL-order.

We no return to the general exposition.
To construct our matching, when we come to a chain
$C:\ w=v_0 \edg v_1 \edg \ldots \edg u$ in a given order
we must be able to determine which faces of $C$ are new.
Denote the open interval $I$ from $v_i$ to $v_j$ in $C$ by
$$
I=C(v_i,v_j)=v_{i+1}\edg v_{i+2}\edg \ldots\edg v_{j-1}.
$$
(Do not to confuse this with an open interval in the poset.)
Then $I$ is a {\it skipped interval\/} if $C-I\sbs C'$ for some
$C'<C$.  It is a {\it minimal skipped interval\/} or {\it MSI\/} if it
does not strictly contain another skipped interval.
In Figure~\ref{322}, the MSI's are circled.  One can find
the MSI's by taking the maximal intervals in $C-(C\cap C')$ for
each $C'<C$ and then throwing out any that  are not containment
minimal in $C$.  Let $\cI=\cI(C)$ be the set of MSI's in $C$.  Then it
is easy to check that a face $\si$ is new in $C$ if and only if
$\si$ has a nonempty intersection with every $I\in\cI(C)$.

The set $\cI$ in not quite sufficient to construct the matching
because the MSI's can overlap and we will need disjoint intervals.
However, there are no containments among the intervals in $\cI$, so
they can be ordered $I_1,I_2,I_3,\ldots$ according to when they are  
first
encountered on $C$.  We now inductively construct a set $\cJ=\cJ(C)$ of
{\it  $J$-intervals\/} as follows.  Let $J_1=I_1$.  Then consider
the intervals $I_2'=I_2-J_1,I_3'=I_3-J_1,\ldots$; throw out any which  
are not
minimal; and pick the first one which remains to be $J_2$.  Continue  
this
process until there are no nonempty modified MSI's left.  It
happens that in all the critical chains for $\bbP^*$, the
intervals in $\cI$ will already be disjoint and so we will not need
this step.

We are finally in a position to describe the matching.
List the maximal chains of $[u,w]$ using a PL-order.
A family $\cK$ of intervals of a maximal chain $C$ {\it covers\/} the
chain if $C=\cup \cK$.
There are three cases depending on whether $\cI(C)$ or $\cJ(C)$ covers
$C$ or not.
First suppose that $\cI(C)$ does not cover $C$, so neither does  
$\cJ(C)$,
and pick $x_0$ to be the first vertex in $C-\cup \cI$.  Consider the map
$\si\ra\si\De\{x_0\}$ where $\De$ is symmetric difference (not the
order complex).
One can show that this map is a fixed-point free involution on
the new faces $\si$ in $C$ which extends the Morse matching already
constructed from the previous chains.
Now suppose that $\cJ=\{J_1,J_2,\ldots,J_r\}$
does cover $C$ and consider the new face $\si_C=\{x_1,x_2,\ldots,x_r\}$  
where
$x_i$ is the first element of $J_i$ for $1\le i\le r$.  Given any
other new face $\si\neq\si_C$, we find the interval $J_i$ of smallest
index where $\si\cap J_i\neq\{x_i\}$ and map $\si\ra\si\De\{x_i\}$.
This involution pairs up all new faces in $C$ except $\si_C$, which is
critical.  Finally, suppose that $\cI$ covers $C$ but $\cJ$ does not.
Then we use the mapping of the second case to pair up all
new faces whose restriction to $\cup\cJ$ is different from $\si_C$.
We also pair up the remaining new faces (including $\si_C$) by using
the mapping of the first case where we take $x_0$ to be the first
vertex in $C-\cup\cJ$.
Thus we have outlined the proof of the following theorem,
remembering that the dimension of a simplex is one less than its
number of vertices.
\bth[\cite{bh:dmf}]
\label{bh}
Let $P$ be a poset and $[u,w]$ be a finite
interval in $P$.  For any PL-order on the maximal chains of $[u,w]$,
the above construction produces a Morse matching in $\De(u,w)$ with
the following properties.
\ben
\item  The maximal chain $C$ is critical if and only if $\cJ$ covers  
$C$.\vs{3pt}
\item  If $C$ is critical then its unique critical cell has
   dimension $\#\cJ(C)-1$.\qed
\een
\eth

\section{A Morse theory derivation of $\mu$}
\label{pdm}

We are now ready to find the critical cells for the PL-order in
$\bbP^*$ defined previously.  We first need three lemmas
which will prove useful in a number of cases.  Unless otherwise
specified, we always use the notation
\beq
\label{C}
C:\ w=v_0 \edge{l_1} v_1 \edge{l_2} v_2 \edge{l_3}\ldots \edge{l_d}  
v_d= u
\eeq
for labeled maximal chains, or
\beq
\label{Ce}
C:\ w=\eta_{v_0} \edge{l_1} \eta_{v_1} \edge{l_2} \eta_{v_2}
\edge{l_3}\ldots \edge{l_d} \eta_u
\eeq
if we wish to be specific about the embeddings determined by $C$.  We  
also use
$$
l(C)=(l_1,l_2,\ldots,l_d)
$$
for its label sequence.

Take an interval $[u,w]\sbs\bbP^*$ with $|u|=|w|$ and let  
$m_i=w(i)-u(i)$.
Now consider the multiset $M_{uw}=\{\{1^{m_1},2^{m_2},\ldots\}\}$ where
$i^{m_i}$ means that $i$ is repeated $m_i$ times.  Then every
permutation of $M$ is the label sequence for a unique maximal chain in
$[u,w]$ and this accounts for all the chains.  (In fact, $[u,v]$ is
isomorphic to the poset of submultisets of $M_{uw}$.)  We record this  
simple
observation for later reference.
\ble[Same Length Lemma]
\label{sll}
If $|u|=|w|$ then the the label function $l$ gives a bijection between
the maximal chains  in $[u,w]$ and the permutations of $M_{uw}$.  In
particular, if $M_{uw}$ contains only one distinct element (possibly
with multiplicity) then $[u,w]$ contains a unique maximal chain.\qed
\ele

If $|u|<|w|$ then we no longer have the nice bijection of the previous
paragraph, but we can still say something.  Let $C$ be a maximal
chain as in~\eqref{C} and let $l'=(l_1',\ldots,l_d')$ be any permutation
of the label sequence $l(C)$.  Then $l'$ defines a sequence of  
expansions
$\eta_{v_0'},\eta_{v_1'},\ldots,\eta_{v_d'}$ where $\eta_{v_0'}=w$ and
for $j\ge1$ we get $\eta_{v_j'}$ from $\eta_{v_{j-1}'}$ by subtracting
one from position $l_j'$ in $\eta_{v_{j-1}'}$.  It is still true that
$C':w=v_0'\edg v_1'\edg\ldots\edg v_d'=u$ is a maximal chain in
$[u,w]$.  We call $C'$ the {\it chain specified by $l'$}.
Since $\eta_{v_j'}$ may not be a normal embedding in
$\eta_{v_{j-1}'}$, we may not have $l(C')=l'$.  Still, at the first
place where $l(C')$ and $l'$ differ, that difference must have been
caused because using the label in $l'$ would have resulted in changing
a $1$ to a $0$ where that $1$ was not the first in its run.  Thus the
corresponding normal embedding in $l(C')$ uses the first $1$ in that
run which is to the left.  Hence $l(C')\le l'$ in lexicographic order.
We summarize this discussion in the following lemma.
\ble[Chain Specification Lemma]
\label{csl}
If $C$ is a maximal chain in $[u,w]$ and $l'$ is any permutation of
$l(C)$ then $l(C')\le l'$ where $C'$ is the chain specified by $l'$.\qed
\ele

As our first application of the Chain Specification Lemma, we can
determine what happens at descents.  A {\it descent\/} of $C$ is
$v_j\in C$ such that $l_j > l_{j+1}$.  An {\it ascent\/} is defined by
reversing the inequality.
\ble[Descent Lemma]
\label{dl}
If $v_j$ is a descent of $C$ then it is an MSI.
\ele
\begin{proof}
Let $l'$ be the permutation of $l(C)$ gotten by interchanging $l_j$
and $l_{j+1}$ and let $C'$ be the chain specified by $l'$.  Then by
Lemma~\ref{csl} we have $l(C')\le l'<l(C)$, so $C'$ comes before $C$
in PL-order and it is easy to check that $C'$ diverges from $C$ at
$v_{j-1}$ and rejoins $C$ at $v_{j+1}$.  Thus $\{v_j\}$ is a skipped
interval; and since the interval contains only one element it must
also be minimal.
\end{proof}

We only need a few more definitions to state our result characterizing
the critical chains.  A chain $C$ will be said to have a certain  
property,
e.g., weakly decreasing, if $l(C)$ has that property.  Also, if
$\eta_u$ is a normal embedding into $w$ then we need to keep track of
the zero positions which did not come from decreasing a $1$ in $w$ by  
letting
$$
D(\eta_u)=\#\{i\ |\ \eta_u(i)=0,\ w(i)\ge2\}.
$$
\bth
\label{morse}
Consider the maximal chains  in $[u,w]\sbs\bbP^*$ in the given
PL-order.
\ben
\item  There is a bijection between critical chains $C$ and normal
   embeddings $\eta_u$ into $w$ where the chain corresponding to $\eta_u$
   is the unique weakly decreasing chain ending at $\eta_u$.
\item  If $C$ is critical and ends at $\eta_u$ then $\cI(C)=\cJ(C)$
   and
$$
\#\cI(C)=d(\eta_u)+2D(\eta_u)-1.
$$
\een
\eth

We shall prove this theorem by considering 3 cases: when $C$ is weakly
decreasing and ends at a normal embedding, when $C$ is weakly
decreasing and does not end at a normal embedding, and when $C$ is not
weakly decreasing.  Note that for any embedding $\eta_u$ into $w$,
there is at most one decreasing chain ending at $\eta_u$, and that if
$\eta_u$ is normal then such a chain will exist because it will be
possible to make each cover normal.  Thus there is a bijection between
normal embeddings and weakly decreasing chains ending at them, but we
need to show such chains are critical.  To do so, we define a
{\it plateau\/} of $C$ to be an interval $C(v_i,v_j)$ such that
$(l_{i+1},l_{i+2},\ldots,l_j)$ is a run of length at least 2 in $l(C)$.

\bpr
\label{critical}
If $C$ is weakly decreasing and ends at a normal embedding $\eta_u$
then $C$ is critical, $\cI(C)=\cJ(C)$, and  
$\#\cI(C)=d(\eta_u)+2D(\eta_u)-1$.
\epr
\begin{proof}
Every descent is an MSI of $C$ by the Descent Lemma, so any other
MSI must be contained in a plateau by minimality.
In fact, we claim that any plateau $C(v_i,v_j)$ is an MSI.
Without loss of generality we can assume $v_i=w$ (since otherwise
$v_i$ is a descent and so no MSI can contain it) and $v_j=u$.

To show $C=C(w,u)$ is a skipped interval, first note that by  
construction
$l(C)$ consists of $c$ repeated $k=j-i\ge2$ times, so  
$w(c)=\eta_u(c)+k$.
Thus by normality and the fact that $k\ge2$ we get
$\eta_u(c)=0$ and $w(c)=k$.  Using  normality again implies that $c$ cannot  
be the
first element in its run of $k$'s in $w$,
and thus $w(c-1)=k$.  Because of this, there is a chain
$C'$ from $w$ to $u$ all of whose labels are $c-1$.
By construction $C'<C$ and $C\cap C'=\emp$, so $C$ is a MSI as desired.

To show the plateau is minimal suppose, to the contrary, that there
is a skipped interval $I\sbs C$ and let $w=v_0$, $u=v_j$.  Note that
because $|v_0|=|v_1|=\ldots=|v_{j-1}|$, the Same Length Lemma applies
to show that there is only one chain (namely an interval of $C$)
between any two of these compositions.  Thus the chain $C'$ giving rise
to $I$ must rejoin $C$ at an embedding $\eta_u'$ of $u$ into $w$, and
hence also contain $v_1$ in order to cut out a proper subinterval.
 From this and normality of $\eta_u$ we have
\begin{equation}\label{thedisplayedequationabove}
\eta_u'(c)<k=w(c-1)=\eta_u(c-1).
\end{equation}
Also, $|u|=|w|-1$ implies that $C'$ must zero out exactly one element
of $w$.  Since $C'<C$, that element must be in a position strictly to  
the
left of position $c$.   But then because $\eta_u$ and $\eta_u'$ are both
expansions of $u$ we are forced to have $\eta_u(c-1)=\eta_u'(c)$,
contradicting \eqref{thedisplayedequationabove}.

Now we know that $\cI(C)$ consists of the descents and plateaus of $C$
which are disjoint and cover $C$ by their definition, so
$\cJ(C)=\cI(C)$ and $C$ is critical by Theorem~\ref{bh}.  To count
the number of MSI's, note that if $c$ is a position counted by
$d(\eta_u)$ then the vertex just before the edge labeled $c$  in $C$
will be a descent, unless that edge is the very first one.  On the
other hand, if $c$ is counted by $D(\eta_u)$ then the run of $c$'s in
$l(C)$ contribute both a plateau and a descent just before the plateau
to $\cI$ (unless the run is at the beginning of $C$ when only the
plateau will be an interval).  In this manner we count
each MSI exactly once for a total of $d(\eta_u)+2D(\eta_u)-1$ intervals.
\end{proof}

\bpr
If $C$ is weakly decreasing and ends at an embedding $\eta_u$ which is
not normal then $C$ is not critical.
\epr
\begin{proof}
As in the proof of the previous proposition, it suffices to consider
the case where $l(C)$ consists of a label $c$ repeated $k\ge2$ times
so that $w(c)=\eta_u(c)+k$.  If $\eta_u(c)>0$ then $|w|=|u|$ and so
the Same Length Lemma applies to show that $C$ is the only chain from
$w$ to $c$.  In particular, it is the lexicographically first chain
and thus not critical.

If $\eta_u(c)=0$ then $w(c)=k$ and, since $\eta_u$ is not normal, it
must be that $c$ is the first index in this run of $k$'s in $w$.  Let
\beq
\label{eta}
\eta_u(c-1)=w(c-1)=h\neq k.
\eeq
To demonstrate that $C$ is not critical, it suffices to show that
there is no MSI $I$ containing the element $v_1$ in $C$.
Suppose, to the contrary that such an interval
$I$ exists and let $C'$ be a chain giving rise to $I$.  Using the
Same Length Lemma as in the proof of Proposition~\ref{critical} (third
paragraph) we see that $C'$ must rejoin $C$ at $u$ and this forces
$I=C(w,u)$.

To finish the proof, it suffices to find a skipped interval
$I''\sbs I$ since that will contradict the minimality of $I$.  Let $b$  
be
the smallest label in $l(C')$.  Then $b<c$ since $l(C')<l(C)=(c^k)$.
Now using the same argument as at the end of the third paragraph of
the previous proposition as well as equation~\eqref{eta} gives
$$
\eta_u'(c)=\eta_u(c-1)=h\neq k=w(c).
$$
Since the parts of a composition can only (weakly) decrease along a
chain, we must have $\eta_u'(c)<w(c)$.  It follows that $c$ is a label
on $C'$.  Now consider any permutation of $l(C')$ which starts
$l''=(c,b,...)$ and let $C''$ be the chain specified by $l''$. Then by
construction  and the Chain Specification Lemma $l(C'')\le l''<l(C)$.
This implies that $C''<C$ in PL-order and, by
construction again, $C''$ contains $v_1$.
Thus no MSI of $C$ can contain $v_1$, a contradiction.
\end{proof}

Our third and final proposition completes the proof of  
Theorem~\ref{morse}
\bpr
If $C$ is not weakly decreasing then $C$ is not critical.
\epr
\begin{proof}
If $C$ is not weakly decreasing then it has an ascent $v$.  It
suffices to show that $v$ is in no MSI.  Suppose, to the contrary,
that $v$ is in an MSI $I=C(v_i,v_j)$.  Then by the Descent Lemma, $I$  
contains no
descents and so $C$ is weakly increasing from $v_i$ to $v_j$.  But as
in the previous two proofs, it is no loss of generality to assume
$v_i=w$ and $v_j=u$ so that $C$ is itself an MSI.

Since $C$ is an MSI, it is not the first chain in $[u,w]$.  That first
chain is the unique weakly increasing chain which ends at the
rightmost embedding $\rho_u$ of $u$ into $w$.  Thus if $\eta_u$ is the
embedding defined by $C$ then we must have $\eta_u\neq\rho_u$.

For $\eta_u$ define
$$
z_\eta(j)=\#\{i\le j\ |\ \eta_u(i)=0\}
$$
and similarly define $z_\rho(j)$ for $\rho_u$.  Because $\rho_u$ is
rightmost we always have $z_\rho(j)\ge z_\eta(j)$ with equality when
$j=|w|$.  But $\rho_u$ is not equal to $\eta_u$, so there is an index
$a$ such that $z_\rho(a)> z_\eta(a)$.  Thus there is a first index
$c>a$ such that $z_\rho(c)= z_\eta(c)$.  This definition of $c$ forces
$\eta_u(c)=0$ and $\rho_u(c)=k$ for some $k>0$.  But $\eta_u$ and
$\rho_u$ are embeddings of the same composition, so there must be some
index $b$ with $a\le b<c$ such that $\eta_u(b)=\rho_u(c)=k$  and
$\eta_u(i)=0$ for $b<i\le c$.

We can now derive a contradiction by constructing a smaller skipped
interval in $C$ as follows.  We have $w(c)\ge\rho_u(c)=k$ and
$\eta_u(c)=0$.  Since $C$ is weakly increasing, the labels equal
to $c$ must occur as a plateau.
Therefore there must be vertices $w',u'\in C$ such that $I'=C(w',u')$  
satisfies
$\eta_{w'}(c)=k$, $\eta_{u'}(c)=0$, and $l(I')=(c^k)$.  But
$$
\eta_{w'}(i)=\eta_{u'}(i)=\eta_u(i)=\case{k}{if $i=b$,}{0}{if $b<i<c$,}
$$
so there is a chain $C'$ from $w'$ to $u'$ with $l(C')=(b^k)$.  Since
$b<c$, $I'$ is a skipped interval and we have obtained the desired  
contradiction.
\end{proof}

We can now rederive the formula for $\mu(u,w)$ in
$\bbP^*$.  Combining equation~\eqref{chit} with
Theorems~\ref{me}, \ref{bh}, and~\ref{morse} we obtain
$$
\mu(u,w)=\chit(u,w)=\sum_C (-1)^{\dim\si_C}
=\sum_{\eta_u}(-1)^{d(\eta_u)+2D(\eta_u)-2}
=\sum_{\eta_u}(-1)^{d(\eta_u)}
$$
where the first sum is over all critical chains $C$ in $(u,w)$ and the
other two are over all normal embeddings $\eta_u$ into $w$.

We end this section by remarking that the Morse method can be used as
a powerful tool not just for proving theorems but for discovering the
correct statement to be proved.  The reader may have found our
definition of a normal embedding somewhat ad hoc.  However,
by starting with the very natural chain labeling used above and
looking at the critical chains, one is quickly led to this definition in
order to characterize the embeddings at which such chains end.
Similarly, the defect may seem to have come out of nowhere, but
in order to
determine the dimension of the critical cells one is forced to define
this quantity as well as its big brother $D(\eta_u)$.

\section{Generalized subword order}
\label{gso}

We can now generalize both our result and Bj\"orner's as follows.  Let
$(P,\le_P)$ be any poset.  {\it Generalized subword order\/} is the  
partial
order on $P^*$ obtained by saying that $u\le_{P^*} w$ if
$w$ contains a subsequence $w(i_1), w(i_2),\ldots,w(i_l)$
such that $u(j)\le_P w(i_j)$ for $1\leq j\leq l=|u|$. 
We get ordinary subword order when $P=A$ is an antichain and we get
the composition poset when $P=\bbP$.

It is a simple matter to recast this generalized order in terms of
embeddings.  Let $\zh$ be a special element which is not in $P$ and
let $\Ph$ be the poset obtained by adjoining $\zh$ as a
minimum element, i.e., $\zh<_{\Ph} x$ for all $x\in P$.  Then the  
definitions
of support and expansion are as
usual, just replacing $0$ with $\zh$.  An {\it embedding\/} of $u$ into
$w$ is a length $|w|$ expansion $\eta_u$ of $u$ with
$$
\mbox{$\eta_u(i)\le_{\Ph} w(i)$ for $1\le i\le|w|$.}
$$
As expected, $u\le_{P^*} w$ if and only if there is an embedding of
$u$ into $w$.

Finding an analogue of normality in this context is more delicate.  So
far, we have only been able to do it for a special class of posets.
But there is evidence that more general results are possible; the
next section contains a discussion of this issue.
First note that the definition of a run carries over verbatim to any  
$P^*$.
Now call $P$ a {\it rooted
tree\/} if its Hasse diagram is a tree with a minimum element.  A
{\it rooted forest\/} is a poset where each connected component of its
Hasse diagram is a rooted tree.    Note that both antichains and
chains are rooted forests.
Note also that if $P$ is a rooted forest
then $\Ph$ is a rooted tree so the following definition makes sense.
If $x\in P$ where $P$ is a rooted forest then let $x^-$ be the element
adjacent to $x$ on the unique path from $x$ to $\zh$ in the Hasse
diagram for $\Ph$.  For a rooted forest, a {\it normal embedding\/} of
$u$ into $w$ is an embedding $\eta_u$ into $w$ satisfying two  
conditions.
\begin{enumerate}
\item  For $1\le i\le |w|$ we have $\eta_u(i)=w(i)$,
$w(i)^-$, or $\zh$.\vs{3pt}
\item  For all $x\in P$ and every run $[r,t]$ of $x$'s
in $w$, we have\vs{3pt}
\begin{enumerate}
\item $(r,t]\sbe \Supp \eta_u$ if  $x$ is minimal in $P$,\vs{3pt}
\item $r\in\Supp \eta_u$ otherwise.
\end{enumerate}
\end{enumerate}

Finally, we nee the definition of {\it defect\/} in this situation,
which is as expected:
$$
d(\eta_u)=\#\{i\ |\ \eta_u(i)=w(i)^-\}
$$
for a normal embedding $\eta_u$ into $w$.  The following theorem is
the promised generalization of Theorems~\ref{bjo}
and~\ref{mobius-main}.  Both
of the two proofs we have given of the special case where $P=\bbP$
generalize easily, with the minimal elements playing the r\^{o}le of
$1$ and the rest functioning like the integers $k\ge2$.
\begin{theorem}\label{mobius-forest}
Let $P$ be a rooted forest.  Then the M\"obius function of $P^\ast$ is
given by
$$
\mu(u,w)=\sum_{\eta_u}(-1)^{\eta_u},
$$
where the sum is over all normal embeddings $\eta_u$ of $u$ into  
$w$.\qed
\end{theorem}

\section{Comments and open problems}
\label{cso}

There are several possible avenues for future research.  We
discuss some of them here.

\subsection{Generating Functions}

As mentioned in the introduction, Bj\"orner and
Reutenauer~\cite{br:rmf} gave another proof of the formula for $\mu$
in $A^*$ using generating functions on monoids.  Let
$\bbZ\spn{\spn{A}}$ denote the algebra of formal series
using the elements in $A$ as noncommutative variables and the integers
as coefficients.  Such a series can be written
$$
f=\sum_{w\in A^*} c_w w
$$
for certain $c_w\in\bbZ$.  For example, given $u\in A^*$ one can  
consider the series
\beq
\label{m}
m(u)=\sum_{w\ge u} {w\choose u}_n w.
\eeq

Bj\"orner and Reutenauer showed that~\eqref{m} is rational for any $u$
and obtained, upon
specialization of the variables,  nice expressions for various
ordinary generating functions associated with the M\"obius function of  
$A^*$.
They also derived  results for the zeta function of $A^*$.
The map $m:A^*\ra \bbZ\spn{\spn{A}}$ can be extended to a
continuous linear endomorphism of $\bbZ\spn{\spn{A}}$.  In fact, the
full incidence algebra of $A^*$ is isomorphic to a subalgebra of
this endomorphism algebra.  Bj\"orner and Reutenauer
give another proof of Theorem~\ref{bjo} using this fact.

It is natural to try and apply these ideas to $\bbP^*$, and more
generally to rooted forests.  This has been done by Bj\"orner and
Sagan~\cite{bs:rmf}.

\begin{figure}
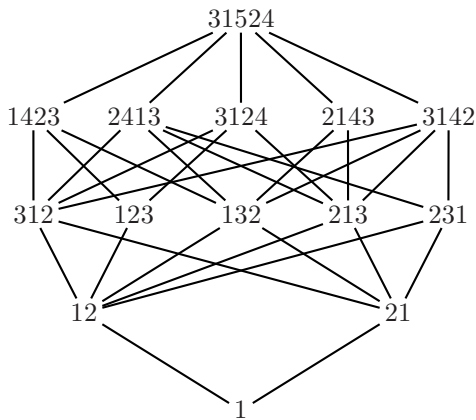

\begin{center}
\begin{footnotesize}
$$
\psset{nodesep=2pt,colsep=18pt,rowsep=25pt}
\begin{psmatrix}
\psspan{2}&[name=31524] 31524&\psspan{2}\\[0pt]
[name=1423] 1423&
[name=2413] 2413&
[name=3124] 3124&
[name=2143] 2143&
[name=3142] 3142\\[0pt]
[name=312] 312&
[name=123] 123&
[name=132] 132&
[name=213] 213&
[name=231] 231\\[0pt]
[name=12] 12\psspan{2}&
&
[name=21] 21\psspan{2}\\[0pt]
\psspan{2}&[name=1] 1&\psspan{2}
\ncline{1423}{31524}
\ncline{2413}{31524}
\ncline{3124}{31524}
\ncline{2143}{31524}
\ncline{3142}{31524}
\ncline{312}{1423}
\ncline{123}{1423}
\ncline{132}{1423}
\ncline{312}{2413}
\ncline{132}{2413}
\ncline{213}{2413}
\ncline{231}{2413}
\ncline{312}{3124}
\ncline{123}{3124}
\ncline{213}{3124}
\ncline{132}{2143}
\ncline{213}{2143}
\ncline{312}{3142}
\ncline{132}{3142}
\ncline{213}{3142}
\ncline{231}{3142}
\ncline{12}{312}
\ncline{12}{123}
\ncline{12}{132}
\ncline{12}{213}
\ncline{12}{231}
\ncline{21}{312}
\ncline{21}{132}
\ncline{21}{213}
\ncline{21}{231}
\ncline{1}{12}
\ncline{1}{21}
\end{psmatrix}
$$
\end{footnotesize}
\end{center}
\caption{The Hasse diagram for the interval $[1,31524]$ in the pattern
containment ordering on permutations}
\label{31524-fig}
\end{figure}

\subsection{The poset of permutations}

Our original interest in $\mathbb{P}^\ast$ came from the
rapidly growing subject of permutation patterns.  For an overview of
permutation patterns the reader is referred to B\'ona's text
\cite{bon:cp}.  Let $S_n$ denote the $n$th symmetric group and let
$\pi\in S_n$ and $\sigma\in S_l$.  We say that 
{\it $\pi$ contains a $\sigma$-pattern\/}, and write $\pi\ge\si$, if
there are indices $i_1<i_2<\dots<i_l$ such that the subsequence
$\pi(i_1)\pi(i_2)\dots\pi(i_l)$ has the same pairwise comparisons as
$\sigma(1)\sigma(2)\dots\sigma(l)$.  This subsequence is called a
{\it copy\/} of $\si$ in $\pi$.  For example, $312\le 24153$
because of  the copy $413$.  This is  a partial
order on the set of all finite permutations.  Wilf was the first to ask
the following question.

\todo{It seems to me that our definition of pattern doesn't need a  
one-line
notation warning, although our example $312\le\dots$ might.}

\begin{question}[Wilf~\cite{wil:pp}]\label{mobius-perms}
What can be said about the M\"obius function of permutations under the
pattern-containment ordering?
\end{question}

Given two permutations $\pi\in S_m$ and $\sigma\in S_n$, their {\it
direct sum\/} is the permutation of length $m+n$ whose first $m$
elements form $\sigma$ and whose last $n$ elements are the
copy of $\pi$ gotten by adding $m$ to each element of $\pi$.
For example, $132\oplus 32145=13265478$.  A permutation is
said to be {\it layered\/} if it can expressed as the direct sum of
some number of decreasing permutations.  (An equivalent
characterization of layered permutations is that they are the
permutations that contain neither a $231$-pattern nor a
$312$-pattern.)  Our previous example is layered because
$13254378=1\oplus 21\oplus 321\oplus 1\oplus 1$.  Clearly the set of
layered permutation of length $n$ is in bijection with the set of
compositions of $n$.  Almost as clearly, this bijection sends the
pattern-containment order to the composition order we have
considered, so Theorem~\ref{mobius-main} answers
Wilf's question for the set of layered permutations.

Any normal embedding approach to describing
the M\"obius function for permutations in general must incorporate
non-unitary weights, as witnessed by the fact that $\mu(1,31524)=6$.

\subsection{Factor order}

Subword order is  not the only partial order on the set of
words.  We say that the word $u$ is a {\it factor\/} of the word $w$ if
there exist (possibly empty) words $v_1$ and $v_2$ so that $w=v_1uv_2$,
or in other words, if $u$ occurs as a contiguous subword in $w$.
Bj\"orner~\cite{bjo:mff} showed that the M\"obius function for factor
order only takes on values in $\{0,\pm 1\}$ and gave a recursive rule
that allows the computation of $\mu(u,w)$ in $O(|w|^2)$ steps.

The factor order can be defined on $P^\ast$ for any poset $P$: we say
that $u$ is a factor of $w$ if there are words $v_1,v_2,v_3$ such that:
\begin{enumerate}
\item $w=v_1v_2v_3$,\vs{3pt}
\item $|v_2|=|u|$,\vs{3pt}
\item $u(i)\le v_2(i)$ for all $1\le i\le |u|$.
\end{enumerate}
Indeed, this is one of the orders on $\mathbb{P}^\ast$ studied by
Snellman~\cite{sne:scc,sne:spa}.  The M\"obius function of
$\mathbb{P}^\ast$ under factor order remains unknown.

\begin{figure}
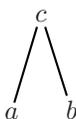

\begin{center}
\begin{footnotesize}
$$
\psset{nodesep=2pt,colsep=18pt,rowsep=25pt}
\begin{psmatrix}
[name=c] c\psspan{2}\\[0pt]
[name=a] a&[name=b] b
\ncline{a}{c}
\ncline{b}{c}
\end{psmatrix}
$$
\end{footnotesize}
\end{center}
\caption{The Hasse diagram for the poset $\Lambda$}
\label{lambda-fig}
\end{figure}
\subsection{Subwords over $\Lambda$}

The smallest poset to which Theorem~\ref{mobius-forest} is
inapplicable is the poset $\Lambda$ depicted in
Figure~\ref{lambda-fig}.  Still, the M\"obius function of $\La^*$
seems to be quite interesting.  In fact, numerical evidence points to a
surprising connection with the Tchebyshev polynomials of the first
kind, $T_n(x)$, which can be
defined as the unique polynomials such that
$$
T_n(\cos\theta)=\cos(n\theta).
$$

\begin{conjecture}\label{mobius-lambda}
For all $i\le j$, $\mu(a^i,c^j)$ is the coefficient of $x^{j-i}$ in
$T_{i+j}(x)$.
\end{conjecture}
As with the poset of permutations, a normal embedding interpretation
of $\mu(a^i,c^i)$ would need to use weights because, for example,
$\mu(a,cc)=-3$.

One possible way to attack this conjecture is to use the three-term
recurrence for $T_n(x)$.  Translating this in terms of the conjecture,
it would suffice to show that
$$
\mu(a^i,c^j)=2\mu(a^i,c^{j-1})-\mu(a^{i-1},c^{j-1})
$$
for $j\ge i\ge1$  However, we have not been able to see any
relationship between the intervals $[a^i,c^j]$, $[a^i,c^{j-1}]$, and
$[a^{i-1},c^{j-1}]$ which would permit us to derive this relation
for their M\"obius functions.

There are two closely related areas where the Tchebyshev polynomials
have appeared.  A permutation $\pi$ {\it avoids\/} a
permutation $\si$ if it doesn't not contain a $\si$-pattern.  Chow and
West~\cite{cw:fsc} showed that the generating function for the number
of elements in $S_n$ avoiding both $132$ and $12\ldots k$ for fixed
$k$ can be expressed in terms of Tchebyshev polynomials of the second
kind.  Mansour and Vainshtein~\cite{mv:rpc} extended this result to
count permutations avoiding $132$ and containing exactly $r$ copies of
$12\ldots k$.

More recently, Hetyei~\cite{het:tp} defined poset maps $T$ and $U$
which he called Tchebyshev transformations of the first and second
kind.  This is because when applied to the ladder poset $L_n$,
the $\bc\bd$-index of the
images can be expressed in terms of $T_n(x)$ and $U_n(x)$.  Since the
$\bc\bd$-index is related to the M\"obius function, it is conceivable
that Hetyei's map could be used to prove our conjecture.   But the
posets $T(L_n)$ are not isomorphic to any of our intervals $[a^i,c^j]$
in general, so it is not clear how to proceed.  However, these maps
are very interesting in their own right and have been further studied
by Ehrenborg and Readdy~\cite{er:ttf}.

\bigskip

\noindent{\it Acknowledgment.}  We are indebted to Patricia Hersh for
useful discussions and references.

\bigskip
\bibliographystyle{acm}
\begin{small}
\bibliography{ref}
\end{small}

\end{document}